\documentclass[12pt, reqno]{amsart}

\usepackage{eucal,color,graphicx}
\usepackage{caption}
\usepackage{subfig}
\usepackage[scientific-notation=true]{siunitx}

\usepackage{tikz,pgfplots}
\usetikzlibrary{arrows.meta}
\usetikzlibrary{backgrounds, calc, through, intersections}
\usetikzlibrary{arrows,chains,matrix,positioning,scopes}
\usetikzlibrary{shapes.arrows}
\usetikzlibrary{shapes.geometric,calc}

\numberwithin{equation}{section}

\theoremstyle{definition}

\theoremstyle{remark}

\newcommand{\og}{\omega}
\newcommand{\om}{\varOmega}

\newcommand{\lm}{\lambda}
\newcommand{\G}{\Gamma}

\newcommand{\rcore}{r_{\text{core}}}
\newcommand{\rclad}{r_{\text{clad}}}
\newcommand{\nclad}{n_{\text{clad}}}
\newcommand{\ncore}{n_{\text{core}}}

\newcommand{\EE}{\vec{\mathcal{E}}}
\newcommand{\EEE}{\vec{{E}}}

\newcommand{\HH}{\vec{\mathcal{H}}}
\newcommand{\HHH}{\vec{{H}}}

\newcommand{\PPP}{\vec{P}}
\newcommand{\Pbg}{\PPP^\text{bg}}
\newcommand{\Pag}{\PPP^\text{ag}}

\newcommand{\ii}{\hat\imath}
\let\d\partial
\renewcommand{\div}{\mathop\text{div}}
\newcommand{\grad}{\ensuremath{\mathop{{\text{grad}}}}}
\newcommand{\curl}{{\ensuremath{\text{curl}\,}}}
\newcommand{\mean}[1]{\langle{{#1}}\rangle}
\newcommand{\slim}{\sum\limits}
\newcommand{\D}{d}

\newcommand{\Isa}{{\mathcal{I}_s}}

\newcommand{\Nt}{N_{\text{total}}}

\newcommand{\Ne}{N_{\text{excited}}}
\newcommand{\Ng}{N_{\text{ground}}}

\newcommand{\dnedt}{\frac{\partial \Ne}{\partial t}}

\newcommand{\sa}{\sigma^{\text{abs}}}
\newcommand{\se}{\sigma^{\text{ems}}}
\newcommand{\ca}{\psi^{\text{abs}}}
\newcommand{\ce}{\psi^{\text{ems}}}

\newcommand{\kR}{\kappa_R}
\newcommand{\gon}{\gamma_1}
\newcommand{\gtw}{\gamma_2}
\newcommand{\gth}{\gamma_3}
\newcommand{\gfo}{\gamma_4}

\renewcommand{\Re}{\text{Re}}
\newcommand{\tnt}{\tilde{N}_{\text{total}}}
\newcommand{\tkr}{\tilde{\kappa}_R}
\newcommand{\cA}{\mathcal{A}}

\newcommand{\veps}{\varepsilon}
\newcommand{\A}[1]{A^{({#1})}}
\newcommand{\K}[1]{K^{({#1})}}

\graphicspath{{./figs/}}

\begin{document}


\title{Simulation of optical fiber amplifier gain using equivalent short
  fibers}

\author{D.~Drake}
\address{Portland State University, PO Box 751, Portland OR 97207,USA }
\email{ddrake@pdx.edu}

\author{J.~Gopalakrishnan}
\address{Portland State University, PO Box 751, Portland OR 97207,USA }
\email{gjay@pdx.edu}

\author{T.~Goswami}
\address{Portland State University, PO Box 751, Portland OR 97207,USA }
\email{tgoswami@pdx.edu}

\author{J.~Grosek}
\address{Directed Energy Directorate, Air Force  Research Laboratory,
  3550 Aberdeen Ave SE, Kirtland Air Force Base,  NM 87117, USA}
\email{jacob.grosek.1@us.af.mil}

\thanks{This work was supported in part by AFOSR grant FA9550-17-1-0090.}

\begin{abstract}
    Electromagnetic wave propagation in optical fiber amplifiers obeys
  Maxwell equations.  Using coupled mode theory, the full Maxwell
  system within an optical fiber amplifier is reduced to a simpler
  model. The simpler model is made more efficient through a new scale
  model, referred to as an equivalent short fiber, which captures some
  of the essential characteristics of a longer fiber.  The equivalent
  short fiber can be viewed as a fiber made using artificial
  (unphysical) material properties that in some sense compensates for
  its reduced length.  The computations can be accelerated by a factor
  approximately equal to the ratio of the original length to the
  reduced length of the equivalent fiber. Computations using models of
  two commercially available fibers -- one doped with ytterbium, and
  the other with thulium -- show the practical utility of the
  concept. Extensive numerical studies are conducted to assess when
  the equivalent short fiber model is useful and when it is not.


\end{abstract}

\maketitle


\section{Introduction}  \label{sec:introduction}

``Scale models'' are ubiquitous in fields such as fluid dynamics.
They are physical or numerical models that preserve some of the
important properties of an object being modeled while not preserving
the original dimensions of the object. The main goal of this paper is
to formulate and study a miniature scale model of an optical fiber
laser amplifier. Our scale model reduces fiber length to increase
computational efficiency. While unable to preserve all properties of
the original electromagnetic solution, our numerical scale model is
able to approximately replicate the original fiber's power
distribution, as we shall see in later sections. After this
introductory section, we will begin by describing a simplified model
of beam propagation in fibers. This model will then be used to derive,
justify, and verify the scale model.

The importance of fiber amplifiers in enabling our current world of
long-distance fiber optics and submarine telecommunications cannot be
overlooked~\cite{Chesn18}.  High power fiber amplifiers also have many
other uses, for example, as defensive speed-of-light weapons.  High
output powers have been achieved by solid-state optical fiber laser
amplifiers~\cite{JaureLimpeTunne13}. Numerical modeling of these
optical devices has also been effectively used by
many~\cite{JaureOttoStutz15, NaderDajanMadde13, SmithSmith11, Ward13}. Yet,
simulation of full length fibers remains cumbersome and far from being
routine. This is because of the long simulation times and the large
computational resources required. Simulations using the full Maxwell
system are too expensive since there are millions of wavelengths
within any realistically long fiber.  An an example, consider the full
Maxwell simulation of Raman gain attempted
in~\cite{NagarGrosePetri19}: more than five million degrees of freedom
was needed to simulate an extremely short fiber containing 80
wavelengths (less than 0.0001~m). Although a full Maxwell model of a
realistically long (10~m) fiber can be written out (and we shall do so
in Subsection~\ref{ssec:full-maxwell-model}), its numerical solution
is beyond the reach of today's simulation capabilities.

Therefore, simplified models form the current state of the art. It is
somewhat surprising how unreasonably effective these models have
proved to be, despite the drastic simplifications used in their
derivation.  The state of the art in fiber amplifier simulation
consists of beam propagation methods using coupled mode theory
(CMT). We shall introduce the reader to the CMT model in
Subsection~\ref{ssec:cmt} as a simplification of the full Maxwell
model. To facilitate cross-disciplinary readership, we make an effort
to enunciate the assumptions behind such simplifications. Even though
it is not common in the optics literature to view CMT in the backdrop
of emerging developments in reduced-order models, one may view it as
essentially an example of a {\em physics-based reduced-order
  model}. Indeed, in CMT, the electromagnetic solution is expressed
using a ``reduced basis'' consisting of transverse guided modes of the
fiber that encapsulates the energy-transport mechanism in fibers.

Even the simplified CMT model is computationally too demanding to ably
assist with the important open issues in the subject today.  One of
these issues is what is currently recognized to be the main roadblock
to power scaling of beam combinable fiber amplifiers, namely the
nonlinear transverse mode
instability (TMI).
TMI can be described as a
sudden breakdown in beam quality at high power operation, first
observed experimentally~\cite{EidamWirthJaure11}.  As pointed out in
the review~\cite{JaureLimpeTunne13}, when attempting to design
highly coherent lasers
capable of sustained high (average) powers, a practically uncrossable
limit was encountered due to the TMI.  After intensive speculations on the
cause of TMI, the prevailing theory seems to be that the cause is a
temperature-induced grating.  
We believe that numerical modeling is essential for investigating the TMI, and other nonlinearities that arise inside fiber amplifiers, since experimental evidence is mostly limited to examining the amplifier output, rather than the onset of physical effects that occur inside of the glass fiber along its length.  
The current difficulty in using numerical models is the
excessive simulation times: indeed any numerical technique used 
must be able to
solve for the electromagnetic field within a long fiber a
vast number of times. 
Given the great computational burden of capturing length scales as small as 10 $\mu$m, and time scales as small as 10 $\mu$sec (for the thermal problem), techniques that further accelerate the numerical simulations have the potential to significantly enhance the ability for computer modeling to inform experimental designs and configurations in a timely manner.
It is our intent to contribute such
an acceleration technique by developing the above-mentioned scale
model (in Sections~\ref{sec:equiv}--\ref{sec:results}). Studies of
its application to TMI investigations are postponed to the future.

The models are tested using two commercially available examples of
doped step-index fibers, one with ytterbium (Yb) doping in the fiber
core, and another with thulium (Tm) doping. Both are examples of large
mode area (LMA) fibers which support more than one guided transverse
core mode. LMA fibers are of great interest since they permit 
greater light amplification per unit length and help mitigate the
onset of other detrimental optical nonlinearities.
Unfortunately, they are also more 
susceptible to the TMI, and hence stand most to benefit from advances
in numerical simulation. Our active gain model for these fibers
utilizes
the population dynamics of Yb and Tm ions.  
Active gain in fiber amplifiers appears as a nonlinear coupling term
between the Maxwell systems for the (1) less coherent ``pump'' light
that supplies energy for amplification, and the (2) highly coherent
``signal'' (laser) light. The gain mechanism involves exciting the outer most electron of the dopant (Yb or Tm) by absorbing the pump light, and producing more coherent signal light via stimulated emission that allows the excited electrons to return to their ground state.
We have included 
a simplified, yet very typical,
mathematical formulation for the dopant ion population dynamics in
Section~\ref{sec:Tm+Yb}. A few of the initial results obtained in this
work for the simpler Yb-doped case were announced earlier in the
conference proceedings~\cite{GopalGoswaGrose19}. We begin by deriving
the CMT model next.

\section{The CMT model}
\label{sec:model}

Physics-based reduced-order models are now being used successfully in
various simulation techniques~\cite{SwiscMainiPeher18}. In this
section, we introduce such a model for an optical fiber amplifier
starting from Maxwell equations. We will start from the Maxwell system
and describe the assumptions that lead us to the simplified CMT model
consisting of a system of ordinary differential equations (ODE).

\subsection{The full Maxwell model}
\label{ssec:full-maxwell-model}

Suppose the optical fiber amplifier to be modeled is aligned so that
it is longitudinally centered along the $z$-axis; the transverse
coordinates will be denoted by $x$ and $y$ in the Cartesian coordinate
system. The core region of the fiber,
$\{ (x,y, z): x^2 + y^2 < \rcore^2\}$, is enveloped by a cladding that
extends to radius $\rclad$. The fiber is a step-index fiber, i.e., its
refractive index is a piecewise constant function that takes the value
$\ncore$ in the core region and $\nclad$ in the cladding region.  
There is usually another polymer coating that surrounds this inner cladding (composed of fused silica); however, this second cladding/coating can readily be neglected for this analysis since the laser light is mostly guided in the fiber core region.
We want to model a continuous wave, weakly guided
($\ncore - \nclad \ll 1 $), polarization maintaining, large mode area
(LMA) fiber. 
There are various
arrangements in which this fiber amplifier could be seeded and
pumped. We consider the {\em co-pumped/clad-pumped} configuration, wherein a
highly coherent laser light -- which we shall refer to as the {\em signal}-- is injected into the fiber core area at the beginning of
the fiber ($z=0$).  The {\em pump} light is injected into the fiber at
$z=0,$ and unlike the signal, it enters both core and cladding.

Let $\EE_s, \HH_s$ and $\EE_p, \HH_p$ denote the electric and magnetic fields of the signal and pump light, respectively.  The signal and pump fields are assumed to be time harmonic of frequencies $\og_s$ and $\og_p$ respectively, i.e.,
\[
  \begin{aligned}
    \EE_\ell (x,y,z,t) 
    & = \text{Re} \Big[ \EEE_\ell(x,y,z)e^{-\ii \og_\ell t} \Big],
    &\quad
    \HH_\ell (x,y,z,t) 
    & = \text{Re} \Big[ \HHH_\ell(x,y,z)e^{-\ii \og_\ell t} \Big],
  \end{aligned}
\]
Here and throughout, we use the subscript $\ell \in \{s, p\}$ to
distinguish between signal and pump fields.  Note that the $\EE_\ell$
and $\HH_\ell$ are real valued while $\EEE_\ell$ and $\HHH_\ell$ are
complex valued. The signal field $\EEE_s, \HHH_s$ and the pump field
$\EEE_p, \HHH_p,$ are assumed to independently satisfy Maxwell
equations, but are coupled through
the electric polarization terms of the form
$\PPP_\ell \equiv \PPP_\ell(\EEE_s, \EEE_p)$, $\ell \in \{s, p\}$,
which appear in the following time-harmonic Maxwell system,
\begin{equation}
  \label{eq:Maxwell-1}
  \begin{aligned}
\curl \EEE_\ell  - \ii \og_\ell \mu_0 \HHH_\ell
& = 0,
\\
\curl \HHH_\ell + \ii \og_\ell \epsilon_0 \EEE_\ell 
& = -\ii\og_\ell \PPP_\ell,
\end{aligned}
\hspace{2cm} \ell \in \{ s, p\},
\end{equation}
where $\epsilon_0$ is the electric permittivity and $\mu_0$ is the
vacuum magnetic permeability.

All interactions between the propagation medium and the
electromagnetic field are modeled through electric polarization terms.
The traditional polarization model includes linear susceptibility,
namely the background material interaction $\Pbg_\ell$ given as a
function of the index of refraction of the medium that the light
propagates through. Other examples of polarization terms include those
that account for linear loss,
active laser gain
($\Pag_\ell$), thermal effects, and optical nonlinearities such as Brillouin scattering, Raman scattering, and Kerr effects. Here we
focus on active gain polarization and the linear background
polarization, namely, 
\[
\begin{aligned}
\Pbg_\ell(\EEE_\ell)
& = \epsilon_0 (n^2 -1) \EEE_\ell,
&&&
\Pag_\ell(\EEE_\ell)
& = - \frac{\ii \epsilon_0 c n}{\og_\ell} g_\ell \EEE_\ell,
&& \quad \ell \in \{ s, p\},
\end{aligned}          
\]
where c is the speed of light and $g_\ell$ is the {\em active gain
  term that depends on $\EEE_\ell$ in some nonlinear } fashion.
Examples of $g_\ell$ are given in Section~\ref{sec:Tm+Yb}.  Typical
optical operating frequencies imply that within a fiber of realistic
length there are several millions of wavelengths.  Even if a mesh fine
enough to capture the wave oscillations is used, the pollution
effect~\cite{BabusSaute00} in wave propagation simulations destroys
the accuracy of finite element solutions at the end of millions of
wavelengths. Hence, without further simplifications, the
above-described full Maxwell model is not a feasible simulation tool
for realistic fiber lengths. We proceed to develop a physics-based
reduced model.

\subsection{Coupled mode theory}   \label{ssec:cmt}

Experiments indicate that the vast majority of the laser signal is contained within the guided core modes of the fiber, and, likewise, most of the pump light is within the guided cladding modes. This is the basis of an electric field ansatz that
  CMT uses.
Before giving the ansatz, let us  eliminate
$\HHH_\ell$ from~\eqref{eq:Maxwell-1}, to obtain the second order
equation
\begin{equation} \label{eqn: simplified maxwell}
	\curl \curl \EEE_\ell - \og_\ell^2 \epsilon_0 \mu_0 \EEE_\ell =
	\og_\ell^2 \mu_0 \PPP_\ell
\end{equation}
solely for the electric field. Substituting
\begin{equation} \label{eqn: total polarization}
	\PPP_\ell = \Pbg_\ell + \Pag_\ell = \epsilon_0 (n^2 -1) \EEE_\ell - \frac{\ii \epsilon_0 c n}{\og_\ell} g_\ell \EEE_\ell
\end{equation}
into~\eqref{eqn: simplified maxwell}, using $c = 1/\sqrt{\epsilon_0 \mu_0}$ and simplifying we get,
\begin{equation}
  \label{eq:Maxwell-2}
  \curl \curl \EEE_\ell - k_\ell^2 n^2 \EEE_\ell + \ii k_\ell n g_\ell \EEE_\ell =0,
\end{equation}
where $k_\ell = \og_\ell/c$ is the wavenumber corresponding to the
frequency $\og_\ell$.

Next, we assume that the electric field $\EEE_\ell$ can be expressed
as
\[ \EEE_\ell (x,y,z) = U_\ell(x,y,z) \hat{e}_x,
\]
i.e., it is linearly polarized in a fixed transverse direction, which
is taken above to be the $x$-direction (where $\hat{e}_x$ denotes the
unit vector in the $x$-direction). Furthermore, since $\EEE_\ell$ has
high frequency oscillations along the $z$-direction, its variations
along the transverse directions may be considered negligible. It is 
therefore standard in optics to neglect $\grad \div \EEE_\ell.$ These
assumptions imply that the vector equation~\eqref{eq:Maxwell-2}
becomes the following scalar Helmholtz equation for $U_\ell$,
\begin{equation}
  \label{eqn: scalar helmholtz}
  -\Delta U_\ell - k_\ell^2 n^2 U_\ell +
  \ii k_\ell n g_\ell U_\ell = 0.
\end{equation}
Due to the high wave number $k_\ell$, even this simplified scalar
field problem is
computationally intensive. We now proceed to further reduce this
scalar model using CMT.

CMT is usually useful in the analysis of the interaction between
several near-resonance guided modes. For step-index fiber waveguides
these modes are called linearly polarized transverse guided core
modes~\cite{Agraw13}, often referred to simply as {\em LP
  modes}. Mathematically speaking, these modes are finitely many
non-trivial functions $\varphi_m(x,y)$, $m=1,2,\ldots, M_\ell$, 
that decay exponentially at the edge of the
cladding region and satisfy
\begin{equation} \label{eqn: eigen problem}
	(\Delta_{xy} + k_{{\ell}}^2 n^2) \varphi_m = \beta_m^2 \varphi_m,
        \qquad
         m = 1, \ldots, M_\ell,
\end{equation}
where $\beta_m$ is the corresponding \textit{propagation constant} and
$\Delta_{xy} = \partial_{xx} + \partial_{yy}$ denotes the transverse
Laplacian {operator}. 
The CMT approach to solve~\eqref{eqn: scalar
  helmholtz} expresses the solution  using  the ansatz
\begin{equation} \label{eqn: cmt ansatz}
	U_{{\ell}}(x,y,z) = \slim_{m=1}^{M_\ell} \A{\ell}_m(z) \varphi_m(x, y) e^{\ii \beta_m z},
\end{equation}
where $\A{\ell}_m(z)$ denotes the complex field amplitude {of mode $m$}. 
{Therefore, the wavenumber ($k_\ell  n$) for the entire electric field envelop ($U_\ell$) is now decomposed into individual propagation constants ($\beta_m$) corresponding to each guided mode, and the field envelop is now decomposed into parts of amplitudes~$\A{\ell}_m$ having transverse profiles  described by $\varphi_m$.}  

Knowledge of the
form of the solution is thus incorporated {\it a priori} into the
ansatz.  In particular, the physical intuition that the
$\varphi_m$-component should oscillate longitudinally at an
approximate frequency of $\beta_m$ is built in. This justifies the
next assumption that $\A{\ell}_m(z)$ is a slowly varying function of $z$
(having built the fast variations in $z$ into the $ e^{\ii \beta_m z}$
term).  Accordingly, for each $\A{\ell}_m$, we neglect the second-order
derivative $d^2 \A{\ell}_m/d z^2$ for all $m=1, \ldots, M_\ell$. Doing so after
substituting~\eqref{eqn: cmt ansatz} into~\eqref{eqn: scalar
  helmholtz} and using~\eqref{eqn: eigen problem} we obtain
\begin{align} \label{eqn: BPM}
	\slim_{m=1}^{M_\ell} \frac{\D \A{\ell}_m}{\D z} \varphi_m \beta_m e^{\ii \beta_m z} 
	& = 
	\frac{1}{2} \slim_{m=1}^{M_\ell} \A{\ell}_m k_{{\ell}} \varphi_m n g_{{\ell}} e^{\ii \beta_m z},
	&& 0 < z < L.
\end{align}

The next step is to multiply both sides of~\eqref{eqn: BPM} by the
complex conjugate of $\varphi_l$, namely $\overline \varphi_l$, and
integrate. We integrate over $\om_z,$ which represents the fiber cross
section having the constant longitudinal coordinate value of~$z$.
Then,  simplifying using
the $L^2(\om_z)$-orthogonality of the modes, 
\begin{align} \label{cmt}
	\frac{\D \A{\ell}_l }{\D z} 
	& = \sum_{m=1}^{M_\ell} e^{\ii (\beta_m - \beta_l)z } \,
   \K{\ell}_{lm}({I_\ell}, I_{{\ell^{c}}})  \; \A{\ell}_m,
	&& 0 < z < L,
\end{align}
for $l = 1, \ldots, M_\ell$,
where $\K{\ell}_{lm}$ is the mode coupling coefficient, given by 
\begin{equation} \label{eq:Klm}
  \K{\ell}_{lm}({I_\ell}, I_{{\ell^{c}}}) = \frac {k_{{\ell}}}{2\beta_l} \int_{\om_z}
	g_{{\ell}}({I_\ell, I_{\ell^{c}}})\,
	n(x,y) \varphi_m(x,y) \overline{ \varphi_l(x,y) } \; dx \,dy,
\end{equation}
{$\ell^{c} \in \{ s,p \} \setminus \{ \ell \}$}, and $I_s, I_p$ denote the signal and pump irradiance,
respectively, 
{which are formulated later in this subsection}.

For the pump light, the number of guided cladding modes is exceedingly
large: $M_p > 10^5$.  Rather than modeling each of these modes,
it is
sufficient to approximate the pump field as a plane wave, which
effectively acts as the composition of all of the pump guided
modes~\cite{NaderDajanMadde13,SmithSmith11}.
Accordingly, we set 
$M_p = 1$ and the normalized mode $\varphi_1^p = (\sqrt{\pi}\rclad)^{-1}$ (without a transverse
dependence). Since the cladding region is many times larger than the
core, the corresponding propagation constant is
estimated as if this mode travels in a uniform medium of refractive
index $\nclad$, i.e., we set 
$\beta_1 = k_p \nclad = \omega_p \nclad /
c$. Then~\eqref{cmt} yields
\begin{equation} \label{pumpODE}
  \frac{\D \A{p}_1}{\D z} = K_{11}^{p}(I_p, I_s)\A{p}_1,
\end{equation}
for $0 < z < L$, where 
\begin{equation} \label{eq:pumpKlm}
  K_{11}^p(I_p, I_s) = \frac{1}{2 \pi \rclad^2 \nclad} \int_{\om_z}
  g_{p}(I_p, I_s) \,\nclad \ dx \,dy = \frac 1 2 \mean{g_p}.
\end{equation}
Here $\mean{g_p} = |\om_z|^{-1} \int_{\om_z} g_p\; dx \,dy$ denotes
the mean value of $g_p$ taken over $\om_z,$ the area of the fiber
cross section out to $r = \rclad$.

The irradiance is proportional to the square of the field envelop
magnitude, $I_\ell = n |U_\ell|^{2} / (\mu_0 c)$. Using~\eqref{eqn:
  cmt ansatz},
\begin{equation}\label{eq:Irradiance}
	{I_\ell(x,y,z) = \dfrac{n}{\mu_0 c} \left|\sum_{m=1}^{M_\ell} \A{\ell}_m(z) e^{\ii \beta_m z} \varphi_m(x,y) \right|^2.}
\end{equation}
{For the pump plane wave, this reduces to}
\[
	{I_p(z) = \dfrac{n}{\mu_0 c \pi \rclad^2} \left| \A{p}_1(z) \right|^2.}
\]
Using the equation~\eqref{pumpODE} and its complex conjugate,
elementary simplifications lead to the following
governing ODE for the pump irradiance:
\begin{equation}\label{eqn:Ip}
	\frac{\D I_p}{\D z} = \mean{g_p} I_p,
\end{equation}
In view of~\eqref{eqn:Ip}, instead of $\A{p}_1$, we shall use $I_p(z)$ as our
pump unknown. There is no need for the amplitude $\A{p}_1$ in the
remainder of the model. Hence from now on, we write $A_m$ for $\A{s}_m$
dropping the superscript. We shall also simply write $M$ for $M_s$
and $K_{lm}$ for $\K{s}_{lm}$.

Next, consider the signal irradiance, namely the $\ell=s$ case
in~\eqref{eq:Irradiance}. To highlight the dependence of $I_s$ on $A_m
\equiv \A{s}_m$, we use $A \equiv [A_1(z), \ldots, A_M(z)]^t$
to collectively denote the set of all signal mode amplitudes
and write
\begin{equation}
  \label{eq:Is}
  I_s \equiv I_s(x,y,z, A) = \dfrac{n}{\mu_0c}
  \left|\sum_{m=1}^M A_m(z) e^{\ii \beta_m z} \varphi_m(x,y) \right|^2.
\end{equation}
Note that
the modes $\varphi_l(x,y)$ and the propagation constants $\beta_l$ may
be precomputed (and the cost of this precomputation corresponds to the
``off-line'' computational cost in this reduced-order model).

{In order to} complete the CMT model (assuming we have expressions for
$g_\ell$), we need to provide initial conditions at $z=0$,
the {beginning of the} fiber. What is usually known is the power
contained in the pump and signal light. 
The initial pump irradiance $I_p^0 = I_p(0)$ can be calculated
from the initial pump power $P_p^0$ provided at the inlet in a
co-pumped configuration, by $I_p^0 = |\om_0|^{-1} P_p^0$.  We assume
that we also
know how the signal light is split into various modes at the inlet,
i.e., we may set 
$A(0)$ to some given $A^0 = [A_1^0, \ldots, A_M^0]^t.$ In
practice, most of the signal power is usually carried in the first
fundamental mode. 

\begin{center}
  \framebox{\parbox{\textwidth}{ %
      To summarize, the CMT model computes
      \[
        Y(z) = [I_p(z), A_1(z), A_2(z), \ldots, A_M(z)]^t, \qquad 0< z<L,
      \]
      {where each $A_m(z)$ is a signal mode amplitude in the fiber core, and $Y(z)$ satisfies}
	 the ODE system
      \begin{subequations}
        \label{eq:summary}
        \begin{align}
          \label{eq:summary-ODE}
          \frac{\D Y}{\D z}
          & =
            \begin{bmatrix}
              \langle g_p(Y)\rangle & 0
              \\
              0 & \phi(z)\cdot K(Y)
            \end{bmatrix}
                  Y,
                                    && 0< z < L,
          \\
          \label{eq:summary-IC}
          Y(0) &  =
                 [I_p^0, A^0]^t
                                    && z = 0,
        \end{align}
      \end{subequations}
      where $ \phi(z)$ is an $M \times M$ matrix defined by
      $\phi_{lm}(z) = e^{\ii (\beta_m - \beta_l) z},$ $ K(Y) $ is a
      matrix of the same size whose $(l,m)$th entry is
      $\K{s}_{lm}({I_s}, I_p)$ defined in~\eqref{eq:Klm}, and
      $\phi(z)\cdot K(Y)$ denotes the Hadamard product of $\phi$ and
      $K$, i.e., $ [\phi \cdot K]_{lm} = \phi_{lm}K_{lm}$.  }}
\end{center}

\section{Thulium and ytterbium doped fiber amplifiers}
\label{sec:Tm+Yb}

Thulium (Tm)-doped fiber amplifiers~\cite{Jacks04, JacksKing99} can
operate in eye-safe laser wavelengths (larger than 1.4~{$\mu$m}) and can
reach an atmospheric transmission window (2.1--2.2~{$\mu$m}). There are
efficient high-power LEDs that operate in the range of
0.79-0.793~{$\mu$m}, which is a peak absorption bandwidth for Tm-doped
fibers. Cross-relaxations and upconversions occur in Tm-doped
amplifiers.  Even though Tm-doped fibers usually have better TMI suppression
compared to other rare-earth ion doped fibers~\cite{SmithSmith11},
ytterbium (Yb)-doped fiber amplifiers have also emerged as excellent
candidates for high power operation due to their high-efficiencies and
low amplified spontaneous emission 
gain. Yb-doped amplifiers are usually pumped at 976~nm and can lase
around 1064~nm very efficiently.  The dynamics of both these ion
populations are explained below. They complete our model by giving
expressions for $g_\ell$ to be used in~\eqref{eq:summary}.

\subsection{Tm-dopant ion dynamics}

\begin{figure}[h]
  \centering
  \begin{tikzpicture}[scale=3, line width=1mm, >=stealth]
    
    \coordinate (3H6) at (0,0);
    \coordinate (3F4) at (0, 0.52);        
    \coordinate (3H5) at (0, 0.826);       
    \coordinate (3H4) at (0, 1.26);        
    \node[left] at (3H6) {$\mathrm{{}^3H_6}$};
    \node[left] at (3F4) {$\mathrm{{}^3F_4}$};
    \node[left] at (3H5) {$\mathrm{{}^3H_5}$};
    \node[left] at (3H4) {$\mathrm{{}^3H_4}$};
    
    \draw (3H6) --++ (1,0) node [right] {ground};
    \draw ($(3H6)+(1.5,0)$) --++ (1,0);
    
    \draw [name path=3f4line]
    (3F4) --++ (1,0)  -- ($(3F4)+(1.5,0)$) --++ (1,0);
    \draw (3H5) --++ (1,0)  -- ($(3H5)+(1.5,0)$) --++ (1,0);
    \draw (3H4) --++ (1,0)  -- ($(3H4)+(1.5,0)$) --++ (1,0);
    
    \draw[blue, dotted, ->, opacity=0.7] ($(3H4)+(2,0)$)--($(3H5)+(2,0)$);
    \draw[blue, dotted, ->, opacity=0.7] ($(3H5)+(2,0)$)--($(3F4)+(2,0)$); 
    \draw[blue,         ->, opacity=0.7] ($(3F4)+(2,0)$)--($(3H6)+(2,0)$)
    node[right, midway] {1700--2100nm};
    
    \draw [purple, ->, opacity=0.7]($(3H6)+(0.75,0)$)--($(3F4)+(0.75,0)$)
    node[right, midway] {1550--1900nm};
    
    \draw[purple, ->, opacity=0.7] ($(3H6)+(0.5,0)$)--($(3H5)+(0.5,0)$)
    node[right, very near end] {1210nm};
    
    \draw[purple, ->, opacity=0.7] ($(3H6)+(0.25,0)$)--($(3H4)+(0.25,0)$)
    node[right, very near end] {793nm};
    
    \draw [red, opacity=0.0, name path=xr]
    ($(3H4)+(0.95,0)$)--($(3H6)+(1.75,0)$);
    
    \fill [red, name intersections={of=xr and 3f4line}]
    (intersection-1) coordinate (intr) circle (0.1pt);
    
    \draw [red, ->, opacity=0.7] ($(3H4)+(0.95,0)$)--(intr);
    \draw [red, ->, opacity=0.7] ($(3H6)+(1.75,0)$)--(intr);
    
    \node at ($(3H6)+(3,0)$) {$N_0$};
    \node at ($(3F4)+(3,0)$) {$N_1$}; 
    \node at ($(3H5)+(3,0)$) {$N_2$};
    \node at ($(3H4)+(3,0)$) {$N_3$};
    
  \end{tikzpicture}
  
  \caption{Simplified diagram of Tm energy levels}
  \label{fig:Tm_manifolds}
\end{figure}
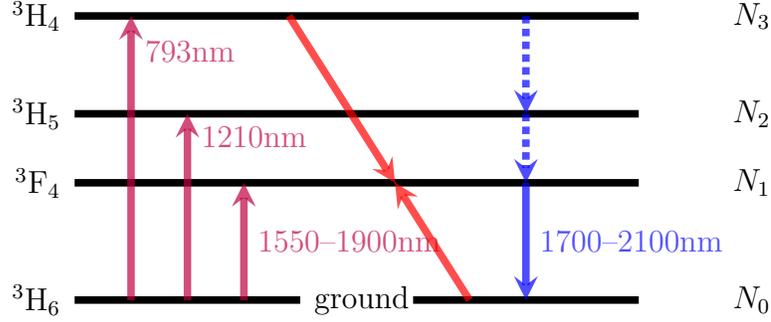

The Tm ion population dynamics are schematically represented in
Figure~\ref{fig:Tm_manifolds}. The model involves four manifolds. The
total number of Tm ions (per volume) is 
\begin{equation}
  \label{eq:Nt-Tm}
\Nt = N_0(x, y, z, t) + N_1(x, y, z, t) + N_2(x, y, z, t) + N_3(x, y, z, t)
\end{equation}
where $N_0$ represents the ground state (manifold 0) ion-population
concentration, while $N_1$, $N_2,$ and $N_3$ denote ion concentrations
at excitation manifolds 1,2 and 3, respectively.  What we have named
energy manifolds 0,1,2, and 3, represent Tm energy levels usually
written as ${}^3H_6$ ${}^3F_4$, ${}^3H_5$ and ${}^3H_4$, respectively.

Pump light of frequency $\og_p=793$~nm excites the Tm ground state
ions into higher energy manifolds, thus depleting manifold 0 at the
rate $\nu_p \sa(\og_p) N_0$ while increasing the excited
manifold $j$ at the rate $\nu_p \se(\og_p) N_j$, where $\sa$ and $\se$
represent measurable absorption and emission cross sections of Tm
\cite{AggerPovls06}, and
\[
  \nu_\ell= \frac{I_\ell} {\hbar  \og_\ell}, \qquad \ell \in \{s, p\}
\]
represents the flux of photons of frequency $\og_\ell$. We must also
take into account the fact that an excited ion in manifold $j$ can
decay spontaneously to a lower energy manifold $k$ at the rate
$1/\tau_{jk}$. An excited ion in manifold $j$ can also decay
non-radiatively to the next lower energy manifold at the rate
$\Gamma_j$. Finally, an excited Tm ion can also undergo
cross-relaxation, wherein it transfers part of its energy to a ground
state ion so both can end up in an intermediate energy
level. Cross-relaxation is represented by the slanted arrows in
Figure~\ref{fig:Tm_manifolds}, while the other processes are
represented by up/down arrows. The rate constant for the
cross-relaxation is denoted by $\kappa_R$.  Cross-relaxation, which
creates two excited Tm ions for every pump photon (a two-for-one
process), increase the amplifier efficiency (while upconversions,
which are neglected in our model, decrease fiber efficiency).
Following~\cite{McCom09}, these processes are modeled by
\begin{subequations}
  \label{eq:Tm-dyn}
  \begin{align}
    \label{eq:Tm-dyn-N3}
    \d_tN_3
    & = \ca_p N_0 - \Big(
      \ce_p  \nu_p + \frac{1}{\tau_{32}} + \frac{1}{\tau_{31}} +
      \frac{1}{\tau_{30}} + \G_3 + \kR N_0 \Big) N_3
    \\     \label{eq:Tm-dyn-N2}
    \d_t N_2
    & = \Big( \frac{1}{\tau_{32}} + \G_3 \Big) N_3 - \Big(
      \frac{1}{\tau_{21}} + \frac{1}{\tau_{20}} + \G_2 \Big) N_2 
    \\     \label{eq:Tm-dyn-N1}
    \d_t N_1
  &= \ca_s N_0 + \Big(
    \frac{1}{\tau_{21}} + \G_2 \Big) N_2 +
    \Big(\frac{1}{\tau_{31}} + 2 \kR N_0 \Big) N_3 -
    \Big( \frac{1}{\tau_{10}} + \G_1 +
    \ce_s \Big) N_1
    \\     \label{eq:Tm-dyn-N0}
    \Nt &  = N_0 + N_1 + N_2 + N_3
\end{align}
\end{subequations}
where
\[
  \ca_\ell = \sa(\og_\ell) \nu_\ell, \qquad \ce_\ell =
  \se(\og_\ell)\nu_\ell, \qquad \ell \in \{s, p\}.
\]
In our simulations, we have set~$\og_s$ to correspond to signal light
of wavelength 2100~nm.

%

Next, we make the simplifying assumption that all the time derivatives
$\d_t$ in~\eqref{eq:Tm-dyn} may be neglected.  By doing so, we are
neglecting the time variations in the ion populations that occur at an
extremely small time scale of around
$10^{-5}$~s. Equations~\eqref{eq:Tm-dyn-N3}--\eqref{eq:Tm-dyn-N1}
after setting $\d_t=0$ immediately yield $N_1, N_2, N_3$ in terms of
$N_0$. The last equation~\eqref{eq:Tm-dyn-N0} then gives a quadratic
equation for $N_0$. To express this solution, first define
\begin{gather*}
  \delta_i    = \slim_{j=0}^{i-1}\tau_{ij} + \G_i,
  \qquad
  \gamma_0 = \frac{1}{\ce_p + \delta_3}
  \qquad 
  \gon   = \ca_p \gamma_0, 
  \qquad 
  \gtw  = \frac{\tau_{32}^{-1} + \Gamma_3}{\delta_2},
  \\
  \gth   = \frac{\tau_{31}^{-1} + \gtw (\tau_{21}^{-1} + \Gamma_2)
    + \gon^{-1} \ca_s}
  {\ce_s + \delta_1},
  \qquad
  \gfo = \frac{2 \ca_p + \ca_s}
  {\ca_p(\ce_s + \delta_1)}.
\end{gather*}
Then, the steady-state solution is given  explicitly by
\begin{subequations}
  \label{eq:N0123}
  \begin{align}
    \nonumber
  N_0 & = \frac{\gamma_0 \kR \Nt - \gon (1 + \gtw + \gth) -1}
        {2 \kR (\gamma_0 + \gon \gfo)}
    \\     \label{eq:N0}
      & +  \frac{\sqrt{(1 - \gamma_0 \kR \Nt + \gon(1 + \gtw +\gth))^2        
    + 4(\gamma_0 + \gon \gfo) \kR \Nt}}{2 \kR (\gamma_0 + \gon \gfo)},
    \\
    \label{eq:N123}
    N_1 &= \frac{(\gth + \gfo \kR N_0)\gon N_0}{1 + \gamma_0 \kR N_0}, \qquad
        N_2 = \frac{\gtw \gon N_0}{1 + \gamma_0 \kR N_0}, \qquad
        N_3 = \frac{\gon N_0}{1 + \gamma_0 \kR N_0}.
\end{align}  
\end{subequations}
Using this, we set the gain expressions by 
\begin{eqnarray}
\label{Tm: signal gain}
g_s &=& \se(\og_s) N_1 - \sa(\og_s) N_0 \\
\label{Tm: pump gain}
g_p &=& \se(\og_p) N_3 - \sa(\og_p)  N_0.
\end{eqnarray} 
This completes the prescription of the CMT model~\eqref{eq:summary}
for Tm-doped fiber amplifier.

\subsection{Yb-dopant ion dynamics}

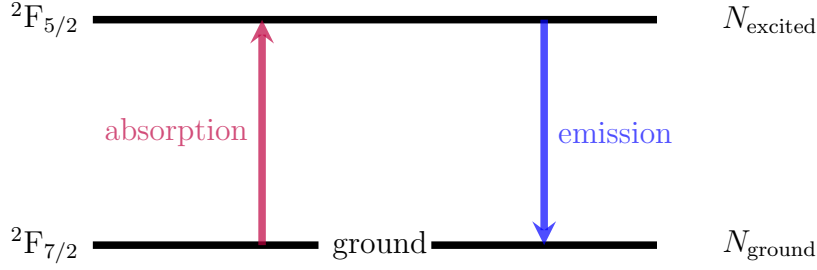
\begin{figure}
  \centering
  \begin{tikzpicture}[scale=3, line width=1mm, >=stealth]
    
    \coordinate (2F72) at (0,0);
    \coordinate (2F52) at (0,1);        
    \node[left] at (2F72) {$\mathrm{{}^2F_{7/2}}$};
    \node[left] at (2F52) {$\mathrm{{}^2F_{5/2}}$};
    \draw (2F72) --++ (1,0) node [right] {ground};
    \draw ($(2F72)+(1.5,0)$) --++ (1,0);
    \draw (2F52) --++ (1,0)  -- ($(2F52)+(1.5,0)$) --++ (1,0);
    
    \draw[blue, ->, opacity=0.7] ($(2F52)+(2,0)$)--($(2F72)+(2,0)$)
    node[right, midway] {emission};
    
    \draw[purple, ->, opacity=0.7]($(2F72)+(0.75,0)$)--($(2F52)+(0.75,0)$)
    node[left, midway] {absorption};
    
    \node at ($(2F72)+(3,0)$) {$\Ng$};
    \node at ($(2F52)+(3,0)$) {$\Ne$};
  \end{tikzpicture}  
  \caption{A simplified diagram of two Yb energy levels}
  \label{fig:Yb_manifolds}
\end{figure}

The model for population dynamics of Yb ions is simpler as it can be
modeled using only two energy states, the ground state and one excited
state manifold, as shown in Figure~\ref{fig:Yb_manifolds}.  Hence,
instead of~\eqref{eq:Nt-Tm}, we now have
$$
\Nt = \Ng(x, y, z, t) + \Ne(x, y, z, t)
$$
where $\Nt$ denotes the total population concentration in the fiber,
$\Ng$ represents the ground state ion-population (in
$\mathrm{{}^2F_{7/2}}$) and $\Ne$ denotes the excited state
ion-population (in $\mathrm{{}^2F_{5/2}}$). The absorption and
emission that models the  two-state dynamics now result in 
\begin{subequations}
  \label{eq:Yb-dyn}
\begin{eqnarray}
  \dnedt & = & \ca_s \Ng - \ce_s \Ne
  \\ \nonumber
  & + &  \ca_p \Ng - \ce_p \Ne - \frac{\Ne}{\tau},
  \\
  \Nt & = & \Ng + \Ne,
\end{eqnarray}   
\end{subequations}
where now we must use the absorption and emission cross section
values~\cite{PaskCarmaHanna95} of Yb for $\sa, \se$ while computing
$\ca_\ell, \ce_\ell$. The parameter $\tau$ is the upper level
radiative lifetime of the excited state. As in the Tm case, we assume
that the system has already reached the steady-state solution. Putting
the time derivative in \eqref{eq:Yb-dyn} to zero, a simple calculation
shows that
\begin{eqnarray}
  \label{steady-state Ne}
  \Ne  & = & \Nt\; \dfrac{ \ca_s +\ca_p}{ \ca_s + \ce_s
             + \ca_p + \ce_p + \tau^{-1}}.
\end{eqnarray}
Finally, the active gain expressions 
are modeled in terms of the above $\Ng$ and $\Ne$  by 
\begin{align}
  \label{gain Yb}
  g_\ell = (\se_\ell \Ne - \sa_\ell \Ng), && \text{for } \ell \in \{ s, p\}.
\end{align}
When this is substituted into~\eqref{eq:summary}, the model for
Yb-doped fiber amplifiers is complete.



\subsection{Basic simulations}   \label{ssec:basic-simulations}

We report the results obtained from simulation of the CMT model for
two 10~m long fibers, one doped with Yb and the other with Tm.  The
fiber parameters are collected from data sheets of commercially
available exemplars of these fibers (specifically
Nufern\texttrademark\ fibers -- see nufern.com).  All parameters used
for the simulation of both the fibers are reported in
Tables~\ref{tab:Yb} and~\ref{tab:Tm}.

\begin{table}
  \centering
  \begin{footnotesize}
  \begin{tabular}{|c|l|l||c|l|l|}
    \hline
    Parameter     & Value           & Units
    & Parameter   & Value           & Units \\
    \hline
    $\lm_p
    =2\pi c/\og_p$  & \num{976e-9}    & m
    & $\lm_s
     =2\pi c/\og_s$ & \num{1064e-9}   & m
    \\
    $\sa(\og_p)$  & \num{1.429E-24} & $\text{m}^2$
    & $\se(\og_p)$     & \num{1.776E-24} & $\text{m}^2$
    \\
    $\sa(\og_s)$  & \num{6E-27}     & $\text{m}^2$
    & $\se(\og_s)$& \num{3.58E-25}  & $\text{m}^2$
    \\
    $\Nt$         & \num{3E+26}     & $\text{ions/m}^3$
    &$\tau$       & \num{8E-4}      & s
    \\
    $\ncore$      & \num{1.450971}  & --
    & NA          & 0.06            & --
    \\
    $\rcore$      & \num{1.25E-5}   & m
    & $\rclad$    & \num{2E-4}      & m
    \\
    $P_p^0$  & 1000           & W
    & $P_s^0$& 25             & W \\
    \hline		
  \end{tabular}
  \end{footnotesize}
  \caption{Parameters used in Yb-doped fiber simulation}
  \label{tab:Yb}  
\end{table}
\begin{table}
  \centering
  \begin{footnotesize}   
    \begin{tabular}{|c|l|l||c|l|l|}
      \hline
      Parameter 		& Value 			& Units 
      & Parameter 		& Value 			& Units\\
      \hline
      $\lm_p=2\pi c/\og_p$	& \num{793E-9} 		& m 
      & $\lm_s=2\pi c/\og_s$ 	& \num{2110E-9} 	& m 
      \\
      $\sa(\og_p)$ 			& \num{4.4686E-25} 	& $\text{m}^2$ 
      & $\se(\og_p)$ 		& 0 				& $\text{m}^2$
      \\
      $\sa(\og_s)$ 			& \num{1.7423E-27} 	& $\text{m}^2$ 
      & $\se(\og_s)$ 		& \num{1.17397E-25} & $\text{m}^2$
      \\
      $\tau_{10}$ 		& \num{6.2232E-03} 	& s 
      & $\tau_{20}$ 	& \num{5.5179E-03} 	& s
      \\
      $\tau_{21}$ 		& \num{2.5707E-01} 	& s 
      & $\tau_{30}$ 	& \num{1.3949E-03} 	& s 
      \\
      $\tau_{31}$ 		& \num{1.7033E-02} 	& s 
      & $\tau_{32}$ 	& \num{6.8446E-02} 	& s 
      \\
      $\G_1$ 			& \num{2.59288E+03} & Hz 
      & $\G_2$ 			& \num{2.92755E+07} & Hz 
      \\
      $\G_3$ 			& \num{8.05943E+04} & Hz 
      & --				& --				& --
      \\
      $\Nt$ 			& \num{3E+26} 		& $\text{ions/m}^3$ 
      & $\kR$ 			& \num{1.17E-21}	& $\text{m}^3$ 
      \\
      $\ncore$ 		& \num{1.439994} 	& -- 
      & NA  		& 0.1 				& --
      \\
      $\rcore$ 		& \num{1.25E-5} 	& m 
      & $\rclad$ 		& \num{2E-4} 		& m 
      \\
      $P_p^0$ 	& 1100 				& W 
      & $P_s^0$ 	& 30 				& W \\
      \hline		
    \end{tabular}
  \end{footnotesize}
    \caption{Parameters used in Tm-doped fiber simulation}
  \label{tab:Tm}  
\end{table}

We solve the CMT system~\eqref{eq:summary} using the classical
$4^{th}$ order explicit Runge-Kutta method (in complex arithmetic).
The phase terms $\phi_{lm}(z) = e^{\ii (\beta_m - \beta_l)z}$ in the ODE
system oscillate at a wavelength not smaller than the so-called {\em
  mode beat length}
\begin{equation}
  \label{eq:mode_beat_len}
  \frac{ 2 \pi}{ \displaystyle \max_{l,m=1, \ldots, M} | \beta_l -
    \beta_m|}.  
\end{equation}
An ODE solver applied to solve~\eqref{eq:summary} must take sufficient
number of steps per mode beat length to capture the effect of these
oscillations in the solution. Prevailing theories
\cite{NaderDajanMadde13} point to the potential importance of the mode
beating term in thermal effects, so we must be careful to treat these
oscillations with the needed accuracy if the model is to be extendable
to incorporate thermal effects in the future.  In all our simulations,
we used 50 ODE steps per mode beat length.

Before running the ODE solver, we precompute the propagation constants
$\beta_j$, the mode beat length, and of course, the modes.  For
step-index fibers, we can compute the modes $\varphi_l$ exactly in
closed form (see~\cite{Agraw13,Reide16}) as quickly described next.
One first computes the propagation constants by solving the
characteristic equation of the fiber as follows.  Let $\mathcal{J}_i$
and $\mathcal{K}_i$ denote, respectively, the standard Bessel function
and the modified Bessel function of second kind of order~$i$. Then we
solve for $X$ satisfying the so-called ``characteristic equation'' of
the fiber, namely setting the fiber's ``numerical aperture''
NA$=\sqrt{\ncore^2 - \nclad^2},$ we solve
$X \mathcal{J}_{i-1}(X) \mathcal{K}_i( \sqrt{\text{NA}^2 - X^2}) +
\sqrt{\text{NA}^2 - X^2} \,\mathcal{J}_i(X) \mathcal{K}_i(
\sqrt{\text{NA}^2 - X^2}) = 0$ by a bisection-based root-finding
method.  This equation arises from the matching conditions at the
core-cladding interface.  For each~$i$, enumerating the roots of the
characteristic equation as $X_{ij}$, $j=0,1,\ldots$, the propagation
constants are given by
\[
  \beta_{ij} = \sqrt{ \ncore^2 k_s^2 - X_{ij}^2 \rcore^2}.
\]
Set $\mathcal{R}_{ij} = X_{ij}/\rcore$ and
$\mathcal{G}_{ij} = \sqrt{\beta_{ij}^2 - \nclad k_s^2 }$.  The exact
LP modes take the following form in polar coordinates:
\begin{equation}
\varphi_{ij}(r, \theta) = 
\begin{cases}
\mathcal{K}_i(\mathcal{G}_{ij} \rcore) \mathcal{J}_i(\mathcal{R}_{ij}
r) \cos (i \theta), & \hspace{1 cm} 0\le r<\rcore
\\
\mathcal{J}_i(\mathcal{R}_{ij} \rcore) \mathcal{K}_i(\mathcal{G}_{ij}
r) \cos (i \theta), & \hspace{1 cm} \rcore \le r < \rclad.
\end{cases}  
\end{equation}
The mode $\varphi_{ij}$ is usually called the ``LP$ij$'' mode.

For the particular case of the Tm parameters in Table~\ref{tab:Tm}, we
find that the fiber only has the LP01 and LP11 modes, while for the Yb
fiber with the parameters set in Table~\ref{tab:Yb}, we found four modes
LP01, LP11, LP21 and LP02. In our simulation the fiber geometry was
meshed using finite elements (with curved elements at the cladding
boundary and at the core-cladding interface) and the relevant LP modes
were interpolated into the degree $p$ Lagrange finite element space
based on the mesh. Integration involving finite element functions is
broken into a sum over integrals over all mesh elements and a
sufficiently high quadrature rule is used to approximated an element
integral. This is how we approximate all required integrals, such as
in the computation of the coupling coefficient~\eqref{eq:Klm}, as well
as in power computations. Note that each step of the multi-stage ODE
solver requires many such integrations.

\begin{figure}
  \centering 
  \includegraphics[width=0.5\textwidth]{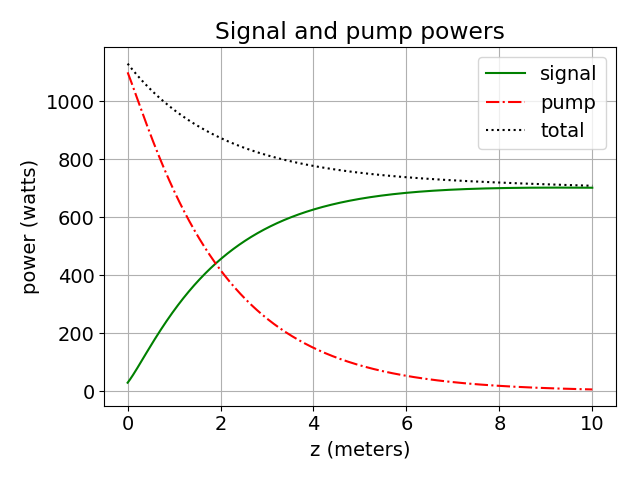}
  \hspace{-0.2cm}
  \includegraphics[width=0.5\textwidth]{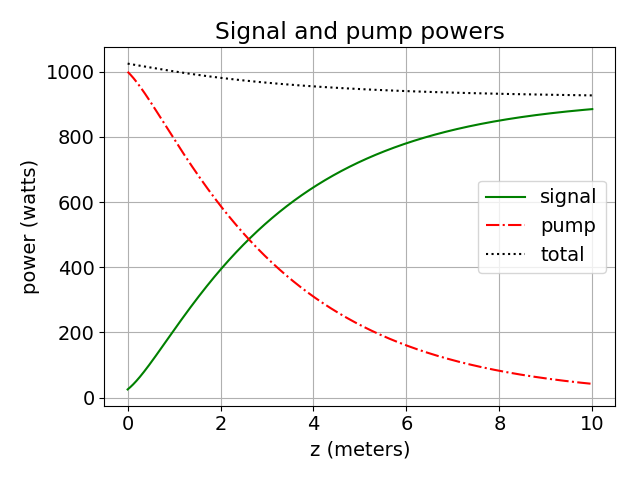}
  \caption{The simulated distribution of powers along the Tm-doped
    (left) and the Yb-doped (right) fiber amplifier. The pump power
    $P_p$ and the signal power $P_s$, as defined in~\eqref{eq:Ps_Pp},
    are shown. The black dotted line plots $P_s+P_p$.}
  \label{fig:Tm_Yb_10m}
\end{figure}

To quantitatively describe the light amplification results of the
simulation, we compute the signal and pump power, after the
approximate $Y(z) = [I_p(z), {A(z)}^t]^t$ has been computed, as follows:
\begin{align}
  \label{eq:Ps_Pp}
  P_s(z)
  & = \int_{\om_z} I_s(x, y, {z}) \;dx dy,
  &&&
      P_p(z)
  & = \int_{\om_z} I_p(z) \;dx dy = |\om_z| I_p(z).
\end{align}          
The initial condition $Y(0)$ is set so that the entire signal power is
fed into the LP01 mode at the inlet $z=0$. Initial pump power $P_p^0$
was set 1000~W for the Yb case and 1100~W for the Tm case.
Figure~\ref{fig:Tm_Yb_10m} shows the distribution of the computed
$P_s$ and $P_p$ (marked ``signal'' and ``pump'' there) for the Tm and
Yb-doped fibers.  The energy transfer from the pump light to the
signal light is clearly evident.  We used $p=5$ Lagrange elements for
these plots.  The use of 50 steps per mode beat length implies that
the Yb case required 421014 RK4 steps, while the Tm case required
302340 steps of the ODE solver to cover the 10~m fiber.

Each of these hundreds of thousands of steps required (multiple)
integrations over the fiber cross section (to compute integrals such
as the one in~\eqref{eq:Klm}). As mentioned above, these integrations
were performed using finite element quadratures. In unreported
experiments, we have attempted to reduce the cost of these
integrations by hyper-reduction techniques common in reduced-order
models~\cite{Rycke09}. One such technique is to use reduced-order
quadratures to approximate the cross-section integrals instead of
using finite elements to perform the integration precisely.  Our pilot
studies into this used Gaussian quadrature rules on a disc (core) and
an annulus (cladding) of order as high as 20.  In cases where this
resulted in substantial reductions in computational cost, we
unfortunately also observed unacceptably large deviations from the
results presented above. Further studies are needed to conclude if
other hyper-reduced quadratures, specifically taking the modes into
account, might prove more useful. In the next section, we describe a
completely different line of inquiry that has yielded considerable
acceleration in our simulations.

\section{The equivalent short fiber concept}
\label{sec:equiv}

In this section, we present the concept of 
{a nearly} equivalent short fiber,
which is an artificially short fiber with unphysical parameters that
can mimic a longer physical fiber in some respects. Being shorter, the
equivalent fiber can be solved using fewer steps of an ODE solver,
thus providing significant reductions in computational cost.

To explain the rationale behind the equivalent short fiber approach,
first consider applying an ODE solver to solve the CMT
model~\eqref{eq:summary}.  As mentioned in the previous section, very
large number of ODE steps were needed to solve the CMT
system~\eqref{eq:summary} on a 10~m long fiber.  Therefore, it would
be extremely useful to reduce the fiber length (and hence the number
of ODE steps) while still preserving the relevant physical processes
in the fiber amplifier. We shall now show that this is possible to
some extent using the computational scale model of an equivalent short
fiber described below.

To begin with, one might consider shortening the $z$-domain
in~\eqref{cmt} using a dimensional analysis.  Note that the left hand
side of~\eqref{cmt} has dimension $\mathrm{V/m}$ (volts per meter),
and $K_{lm}$ has units of $\mathrm{m}^{-1}$. Therefore, by
non-dimensionalization, one is led to believe that a shorter fiber of
length $\tilde L \ll L$ might, in some ways, behave similarly to the
original fiber of length $L$, provided its coupling coefficient is
magnified by $L / \tilde{L}$. However, not all nonlinear systems admit
scale models that are perfect replicas of the original. Below we shall
identify what properties of such a shorter fiber can be expected to be
close to the original.

We introduce the variable change
\[
  \zeta (\tilde z) = \tilde z L/ \tilde L.
\]
A fiber of length $L$, under the variable change
$\tilde z = \zeta^{-1}(z) = z \tilde L/ L$ becomes one of length
$\tilde L$.
Under this variable change, \eqref{cmt} and \eqref{eqn:Ip} become
\begin{eqnarray} 
		\frac{\tilde{L}}{L} \; \frac{\D}{\D \tilde{z}} A_l \Big( \frac{\tilde{z} L}{\tilde{L}} \Big) 
		&=& 
		\slim_{m=1}^M e^{\ii (\beta_m - \beta_l) \tilde{z} L / \tilde{L}} \; K_{lm} \Big( A \Big( \frac{\tilde{z} L}{\tilde{L}} \Big), I_p \Big( \frac{\tilde{z} L}{\tilde{L}} \Big) \Big) \; A_m \Big( \frac{\tilde{z} L}{\tilde{L}} \Big) \\	
		\frac{\tilde{L}}{L} \; \frac{\D}{\D \tilde{z}} I_p \Big( \frac{\tilde{z} L}{\tilde{L}} \Big)
		&=&
		\mean{g_p} \; I_p \Big( \frac{\tilde{z} L}{\tilde{L}} \Big)
\end{eqnarray}
for all $ 0 < \tilde{z} < \tilde{L}$.  In other words, defining
$\hat A_l = A_l \circ \zeta$ and $\hat I_p = I_p \circ \zeta$, the
above system may be rewritten as the following system on the shorter
domain $ 0 < \tilde{z} < \tilde{L}$ for
$\hat Y = [\hat{A}_l, \hat{I_p}]^t$,
\begin{align}    \label{eqn: exact equiv cmt}   
     \frac{\D \hat{Y}}{\D \tilde{z}} 
     &=
       \begin{bmatrix}
       (L/\tilde{L})\; \mean{g_p (\hat Y)} \;
       \hat{I}_p
       \\\displaystyle\slim_{m=1}^M e^{\ii (\beta_m - \beta_l) \tilde{z}
         L / \tilde{L}} \;  (L/\tilde{L}) \; K_{lm} ( \hat Y )
       \; \hat{A}_m, 	
       \end{bmatrix}.
   \end{align}
 Supplemented with the same initial data at at $ z = \tilde{z} = 0$,
\eqref{eqn: exact equiv cmt} is exactly equivalent
to~\eqref{eq:summary}, i.e.,
\begin{equation}
  \label{eq:replica}
  \hat Y = Y \circ \zeta.  
\end{equation}
In other words, the solution of~\eqref{eqn: exact equiv cmt}, being
the pull back of the original solution $Y$ to the shorter domain, is a
perfect replica of the original solution $Y$.

Unfortunately, ~\eqref{eqn: exact equiv cmt} on
$0< \tilde z < \tilde L$ offers no computational advantages over the
original system~\eqref{eq:summary} on $0<z<L$.  This is because the
mode beat length of \eqref{eqn: exact equiv cmt} has been reduced by a
factor of $\tilde{L}/L$ due to the variable change.  So in order
to solve the ODE system \eqref{eqn: exact equiv cmt}, keeping the same
number of steps per mode beat length, the total number of steps needed
to solve the system has not been reduced.
This leads us to consider another mode coupling system with the same
mode beat length as the original system~\eqref{eq:summary}.

\begin{center}
  \framebox{\parbox{\textwidth}{ %
      Let 
      $
      \tilde{Y} (\tilde z) = [\tilde I_p(\tilde z), \tilde{A}_1(\tilde z),
      \cdots, \tilde{A}_M(\tilde z)]^t
      $ solve 
      \begin{subequations}
        \label{eqn:equiv cmt}
        \begin{align}
          \label{eqn:equiv-cmt-1}
          \frac{\D \tilde Y}{\D \tilde z} 
          & =
            \begin{bmatrix}
              \langle (L/\tilde{L}) g_p(\tilde Y)\rangle & 0
              \\
              0 & p(\tilde z)\cdot (L/\tilde{L}) K(\tilde Y)
            \end{bmatrix}
                  \tilde Y,
                                                         && 0< \tilde z < \tilde L,
          \\
          \tilde Y(0) &  =
                 [I_p^0, A^0]^t
                                                         && \tilde z = 0.
        \end{align}  
      \end{subequations}
    }}
\end{center}

Clearly, \eqref{eqn:equiv cmt} is not the same as~\eqref{eqn: exact
  equiv cmt} due to the differences in the phase factors. Therefore,
unlike the solution $\hat Y$ of~\eqref{eqn: exact equiv cmt}, the
solution $\tilde Y$ of~\eqref{eqn:equiv cmt} is not a perfect replica
of the original solution $Y$.  Nonetheless, we shall now proceed to
argue that~\eqref{eqn:equiv cmt} is a practically useful scale model
of~\eqref{eq:summary} as it approximately preserves the power
distribution from the original.  Power, unlike {the amplitude} $A$, is the {quantity that can be, and actually is, experimentally measured}.

Let $P_l$ and $\tilde{P}_l$ be respectively  the powers contained in
the $l^{th}$ mode for the physical and equivalent fiber, defined by 
\begin{eqnarray*}
  P_l(z) 
  & = 
    \displaystyle\int_{\om_z} \frac{n}{\mu_0 c} |A_l (z) \varphi_l(x,y)|^2 \; dx\,dy,
  &
    \qquad \text{ } 0 < z < L,
  \\
  \tilde{P}_l (\tilde z)
  & =
    \displaystyle
    \int_{\om_z} \frac{n}{\mu_0 c} |\tilde{A}_l (z) \varphi_l(x,y)|^2 \; dx\,dy,
  &
    \qquad \text{ } 0 < \tilde z < \tilde L.
\end{eqnarray*}
One may express these in terms of
\[
  \Phi_l = \int_{\om_z} \frac{n}{\mu_0 c} |\varphi_l|^2\;
dx\,dy,
\]
as $ P_l(z) = |a_l|^2 \Phi_l$, where
$ a_l(z) = A_l(z) e^{\ii \beta_l z}$.

To obtain an equation for $P_l(z)$, we may start  from the second equation of
the block system~\eqref{eq:summary}, or equivalently from \eqref{cmt},
which can
be rewritten as 
\[
  e^{\ii \beta_l z}
    d A_l/d z
    = \sum_{m=1}^M K_{lm}(z) e^{\ii \beta_m z} A_m(z).
\]
Then using $da_l/dz = e^{\ii \beta_l z} \d_z A_l + \ii \beta_l a_l,$
we have
\[
  \frac{d a_l}{d z} = \ii \beta_l a_l + \sum_{m=1}^M K_{lm}(z) a_m(z).
\]
Using also the complex conjugate of this equation, 
we have 
\begin{align*}
  \frac{d |a_l|^2}{d z}
  & =   a_l \frac{d \overline{a}_l }{d z}
    + \overline{a}_l \frac{d a_l }{d z}
    =
    \ii \beta_l a_l \overline{a}_l
    - \ii \beta_l  \overline{a}_l a_l
    + \sum_{m=1}^M \overline{K}_{lm} a_l \overline{a}_m +
    {K}_{lm} \overline{a}_l {a}_m,
\end{align*}
i.e.,
\[
  \frac{\D |a_l|^2}{ \D z}
  =
  2\sum_{m =1}^M \Re \big[ K_{lm}(Y) \,\overline{a}_l a_m\big],
\]
for all $l=1, \ldots, M$,  or equivalently,
\begin{equation}
  \label{eqn: pwr cmt}
  \frac{\D P_l}{\D z}
  = 2K_{ll}(Y) P_l  + \rho_l(Y),
\end{equation}
where
\begin{equation}
  \label{eq:rho}
  \rho_l(Y) = 2\Phi_l \sum^M_{\substack{m = 1 \\ m \ne l}}
  \Re \big[ K_{lm}(Y) \,\overline{a}_l a_m\big],
\end{equation}
for $l =1, \ldots, M$.

To the system~\eqref{eqn: pwr cmt}, let us also add the pump power
using the index $l=0$, i.e., let $P_0(z) \equiv P_p(z)$ as defined
in~\eqref{eq:Ps_Pp}.  Then integrating~\eqref{eqn:Ip}, we obtain
$ d P_0/d z = \langle g_p \rangle P_0.$ All together, we have thus
obtained an equation for $P_l$ for all $l=0, \ldots, M$,
\begin{equation}
  \label{eq:Power}
  \frac{d P }{d z}
  =
  \begin{bmatrix}
    \langle g_p(Y)\rangle & 0 
    \\
    0 & 2 \text{diag} [K(Y)]
  \end{bmatrix}
  P 
  +
  \begin{bmatrix}
    0 \\
    \rho(Y)
  \end{bmatrix},
\end{equation}
where $P = [P_0, P_1, \ldots, P_M]^t$ and $\text{diag}[\cdot]$ denotes
the diagonal part of a matrix.

To understand the motivation for the remaining arguments, we now
highlight an observation concerning~\eqref{eq:Power}.  A scale model
providing a perfect replica of the original power distribution is easy
to obtain if the system~\eqref{eq:Power} were an autonomous system:
indeed, if there exists a function $F$ of $P$ alone such that
$d P/ dz = F(P)$, then by merely scaling $F$ by $L/\tilde{L}$, we
obtain an equivalent system that provides perfect replicas of the
original power distribution on the shorter fiber of length $\tilde
L$. However~\eqref{eq:Power} is not autonomous, in general.  Yet, for
practical fibers, our numerical experience suggests that 
\eqref{eq:Power} behaves almost like an autonomous
system. Therefore our strategy now is to view~\eqref{eq:Power} as a
perturbation of an autonomous system.

{Of particular interest is the fact that if the fiber amplifier was robustly single-mode ($M = 1$ for the laser signal), then the governing system}~\eqref{eq:Power} {would be autonomous.  
This can be achieved by not using a LMA amplifier, but one of a smaller fiber core size and/or a lower numerical aperture (NA) such that the fiber core can only support only one guided core mode, the fundamental mode (indexed by $m = 1$), at the signal wavelength.  
However, even with a LMA fiber, if one were to account for fiber bending effects, which cause the higher-order core modes (indexed by $1 < m \leq M$) to leak into the cladding region more so than for the fundamental mode, then the fiber would operate nearly as a single-mode fiber.  
Actual fiber amplifiers are almost always wrapped on a spool rather than stretched out straight, thus ensuring this fiber bending effect.  
This provides us with greater confidence of autonomous system-like behavior, even in real-world implementations of fiber laser amplifier systems.}

Recall from~\eqref{eq:Klm} that $K_{lm}$ is defined using
$ g_s(I_s, I_p)$, where $I_s$ takes the form
in~\eqref{eq:Irradiance}.  We define the following perturbation of $I_s$,
\[
  \Isa(P) =
  \sum_{m=1}^M \frac{n}{\mu_0c}   \left| a_m \varphi_m \right|^2
  =
  \sum_{m=1}^M \frac{n}{\mu_0c \Phi_m}   P_m \left|\varphi_m \right|^2.
\]
It seems difficult to characterize when $I_s - \Isa$ is small {\it a
  priori} (as it depends, e.g., on the localization and orthogonality
of the specific fiber modes) but after a CMT calculation, we may check
if this {difference} is small {\it a posteriori}. Deferring for the moment the
matter of the size of $I_s - \Isa$, let us proceed to
define
$\gamma_\ell(P) = g_\ell(\Isa(P), I_p) = g_\ell(\Isa(P),
P_0/|\om_z|),$ for $\ell \in \{s, p\}.$ They represent the gain
functions obtained by replacing $I_s$ by $\Isa$. The new gain
functions in turn prompt the definition of a new mode coupling
coefficient: instead of~\eqref{eq:Klm}, we now consider
\[
  \kappa_{lm}(P) = \frac {k_s}{2\beta_l} \int_{\om_z}
  \gamma_s(P)
  \,
  n(x,y) \varphi_m(x,y) \overline{ \varphi_l(x,y) } \; dx \,dy.
\]
for all $l, m =1, \ldots, M$. Additionally let 
\[
  \kappa_{00}(P) = \frac 1 2  \mean{\gamma_p(P)} P_0,
\]
and $\kappa_{0l} = \kappa_{l0} = 0,$ for all $l=1, \ldots, M$. We may
now view these $\kappa_{lm}$ as entries of an $(M+1)\times (M+1)$
matrix, using which \eqref{eq:Power} can be expressed as 
\begin{align}
  \label{eq:dPdz}
  \frac{d P }{d z}
  & =
    2\kappa(P) P  + \eta
\end{align}
where   $\eta \in \mathbb{R}^{M+1}$ is defined
by
\[
  \eta(z) =
  \begin{bmatrix}
    \langle g_p(Y) - \gamma_p(P)\rangle  &  0 \\
    0  & 2\,\text{diag}[K(Y) - \kappa(P) ]
  \end{bmatrix}
  P +
  \begin{bmatrix}
    0 \\ \rho(Y) 
  \end{bmatrix}.
\]
We view $\eta$ as a function of $z$, i.e.,
$\eta: [0, L] \to \mathbb{R}^{M+1}$. The $z$-dependence is clear once we
express the $z$-dependence of the solution $Y \equiv Y(z)$ and power
$P \equiv P(z).$ Equation~\eqref{eq:dPdz} shows that power is governed
by a perturbation of an autonomous system whenever $\eta$ is small
enough to be viewed as a perturbation.

Returning to consider~\eqref{eqn:equiv cmt}, we define analogous
quantities for the short fiber, namely
\[
  \tilde{a}_l(z) = \tilde{A}_l(z) e^{\ii \beta_l z}, 
  \quad
  \tilde{P}_0 = \int_{\om_z} \tilde{I}_p \; dx\,dy,
  \quad
  \tilde{P}_l = |\tilde{a}_l |^2 \Phi_l,
\]
for $l=1, \ldots, M$.  Then we may repeat the above arguments starting
from \eqref{eqn:equiv cmt} to obtain the following analogue of~\eqref{eq:dPdz}.
\begin{align}
  \label{eqn: equiv cmt prtrb}
  \frac{\D \tilde{P}}{\D \tilde z}
  & =
    2\frac{L}{\tilde L} \kappa(\tilde P) \tilde P + \tilde \eta,
\end{align}
where $\tilde{\eta} : [0, \tilde{L} ] \to \mathbb{R}^{M+1} $ is now
given by 
\[
  \tilde\eta
    = \begin{bmatrix}
    \langle g_p(\tilde Y) - \gamma_p( \tilde P)\rangle  &  0 \\
    0  & 2\,\text{diag}[K(\tilde Y) - \kappa(\tilde P) ]
  \end{bmatrix}
  \tilde  P +
  \begin{bmatrix}
    0 \\ \rho(\tilde Y) 
  \end{bmatrix}.
\]
Note that $\rho(\tilde Y)$ is defined by~\eqref{eq:rho} after
replacing not only  $Y$ by $\tilde Y,$ but also
$a_l$ (which depends on $Y$) by $\tilde{a}_l$ (which depends on
$\tilde Y$).

To conclude this analysis, it now suffices to compare~\eqref{eqn:
  equiv cmt prtrb} and~\eqref{eq:dPdz}.  Applying the change of
variable $\zeta$ to~\eqref{eq:dPdz}, we get
\begin{equation}
  \label{eqn: equiv pwr cmt var chng}
  \frac{\D }{ \D\tilde z} (P \circ \zeta)
  = 2\frac{L}{\tilde L} \kappa( P \circ \zeta)  P \circ \zeta +
  \frac{L}{\tilde L} \; \eta \circ \zeta.
\end{equation}
Comparing \eqref{eqn: equiv cmt prtrb} and \eqref{eqn: equiv pwr cmt
  var chng} we see that when $\eta$ and $\tilde \eta $ are negligibly
small compared to the other terms, $P_l \circ \zeta$ and $P_l$ solve
approximately the same equation, and consequently
\begin{equation}
  \label{eq:1}
  P\circ \zeta \approx \tilde{P}.
\end{equation}
We summarize this discussion as follows.
\begin{center}
  \framebox{\parbox{\textwidth}{ %
      The system \eqref{eqn:equiv cmt} is an equivalent short fiber
      model of \eqref{cmt} in the sense {that the power $P_l$
        contained in the $l^\text{th}$ mode is approximately preserved
        from the original fiber model \eqref{cmt} through a change of
        variable, under the above assumptions}.  }}
\end{center}



\section{Computational verification of equivalent fiber concept}
\label{sec:results}

In this section, we perform extensive numerical experiments to verify
the pratical utility of the equivalent fiber concept introduced in
Section~\ref{sec:equiv}. We shall compare the relative differences in
the powers obtained from the original fiber and its equivalent short
fiber for various settings to gauge the practical effectiveness of the
approximation~\eqref{eq:1}. In Subsections~\ref{ssec:Tm-equiv-realize}
and~\ref{ssec:Yb-equiv-realize}, we show a way to understand the
equivalent short fiber as a fiber with artificial parameters (with
values not physically realizable) for the Tm and Yb cases,
respectively.

\subsection{Realizing the equivalent short fiber for the Tm-doped  case}
\label{ssec:Tm-equiv-realize}

The equations of the equivalent short fiber, namely~\eqref{eqn:equiv
  cmt}, can be realized for a dopant medium if we can find a set of
``artificial'' parameters that would scale the original $g_p$ and the
original $K$ by $L/\tilde L$. In view of~\eqref{eq:Klm}, this effect
is achieved by scaling the original $g_\ell$ by $L/\tilde L$ for
$\ell \in \{s, p\}$.  Now consider the expressions for $g_\ell$ for
Tm-doped fiber, given in \eqref{Tm: signal gain} and \eqref{Tm: pump
  gain}. Clearly, in view of these expressions, $g_\ell$ will scaled
by $L/\tilde L$ if all the ion populations $N_i$ are so scaled.

This observation, in turn, leads us to consider the expressions for
$N_i$ we derived in~\eqref{eq:N0123}. Let
\[
  \tnt = \frac{L}{\tilde L}  \Nt, \qquad
  \tkr = \frac{\tilde L}{L} \kappa_R.
\]
The value of the expression for $N_0$ in~\eqref{eq:N0}
will be scaled by $L /\tilde L$ if we replace $\kappa_R$ by $\tkr$ and
$\Nt$ by $\tnt$, i.e., \eqref{eq:N0} implies
\begin{align}
  \label{eq:Nt0}
  \frac{L}{\tilde L} N_0
  & = \frac{\gamma_0 \tkr \tnt - \gon (1 + \gtw + \gth) -1}
    {2 \tkr (\gamma_0 + \gon \gfo)}
  \\          \nonumber 
  & +  \frac{\sqrt{(1 - \gamma_0 \tkr \tnt + \gon(1 + \gtw +\gth))^2        
    + 4(\gamma_0 + \gon \gfo) \tkr \tnt}}{2 \tkr (\gamma_0 + \gon \gfo)}.
\end{align}
Let $\tilde{N}_0 = L N_0 / \tilde{N_0},$ the left hand side above.
Proceeding to analyze the expressions in~\eqref{eq:N123}, we find that the
same change in $\kappa_R$ and $\Nt$, and the consequent change in
$N_0$ to $\tilde{N}_0$ per~\eqref{eq:Nt0}, also scales all other $N_i$ by
$L /\tilde L,$ i.e.,
\begin{align*}
  \frac{L}{\tilde L} N_1
  &= \frac{(\gth + \gfo \tkr \tilde{N}_0)\gon \tilde{N}_0}{1 + \gamma_0 \tkr \tilde{N}_0}, \quad
    \frac{L}{\tilde L}  N_2 = \frac{\gtw \gon \tilde{N}_0}{1 + \gamma_0 \tkr \tilde{N}_0}, \quad
    \frac{L}{\tilde L} N_3 = \frac{\gon \tilde{N}_0}
    {1 +  \gamma_0 \tkr \tilde{N}_0}.
\end{align*}
Therefore, all the ion populations $N_i$ are scaled by $L/\tilde{L}$,
and so are $g_s$ and $g_p$.  We have thus arrived at our main
observation of this subsection:
\begin{center}
  \framebox{\parbox{\textwidth}{ %
      A short fiber of length $\tilde L$ is equivalent to a Tm-doped
      fiber of length $L$ if the fiber's original parameters $\Nt$ and
      $\kappa_R$ are changed to $\tnt = L \Nt / \tilde L$ and
      $\tkr = \tilde L \kappa_R /L$, respectively, i.e., this change
      realizes~\eqref{eqn:equiv cmt}.
    }}
\end{center}

\begin{figure}
  \centering
  \includegraphics[width=0.5\textwidth]{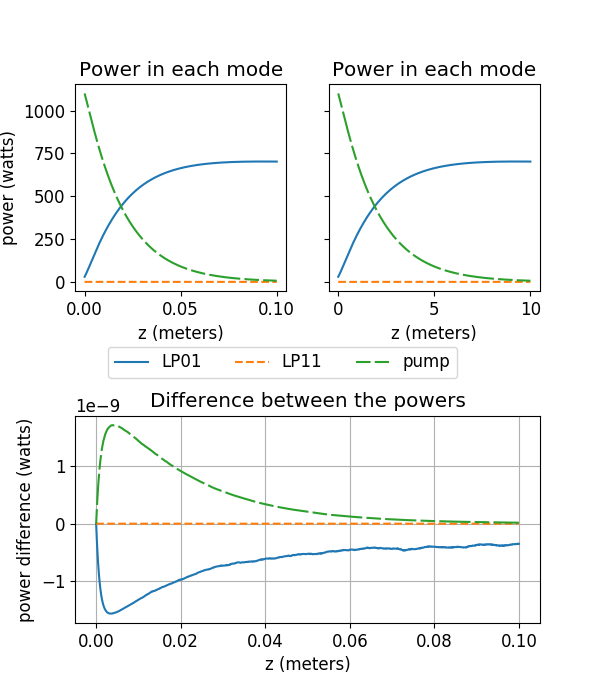}
  \hspace{-0.2cm}
  \includegraphics[width=0.5\textwidth]{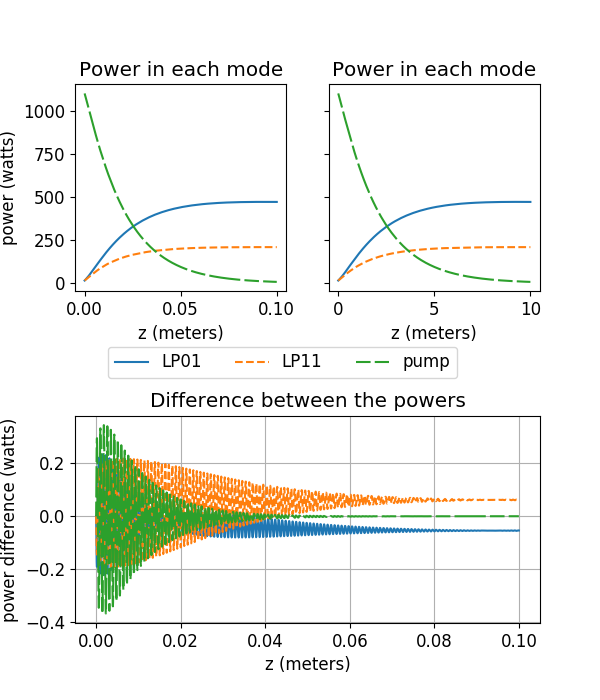}
  \caption{A comparison between a  {\bf{Tm-doped fiber}} and its
    equivalent short counterpart. The left panel shows the case where
    the input signal power was wholly contained in the LP01 mode,
    while the right panel shows the case where it was equally
    distributed between the two modes.}
  \label{fig:Tm-individual}
\end{figure}

To see how this idea works in practice, we consider two scenarios,
both with an equivalent short fiber of $\tilde{L} = 0.1$~m
representing the 10~m long Tm fiber we simulated in
Figure~\ref{fig:Tm_Yb_10m}. (All parameters are as in
Table~\ref{tab:Tm} except for $\Nt$ and $\kappa_R$, which were
modified for the equivalent fiber as stated above.) In the first
scenario, 100\% of the input signal power is carried in the LP01 mode
at the inlet (the same setting as in the computation reported in
Figure~\ref{fig:Tm_Yb_10m}). In the left panel of
Figure~\ref{fig:Tm-individual}, we find that the plots of the computed
powers for the equivalent short fiber and the real fiber are virtually
identical. Even though the difference between them appear to be zero
visually, we have quantified this difference in the bottom left plot
of Figure~\ref{fig:Tm-individual}: since the domains of the two power
functions to be compared are different, we pull back the original
powers to the shorter domain and plot $P_l \circ \zeta - \tilde{P}_l$
(for the two modes, LP01 and LP11) on the shorter domain. Clearly,
from the scale of the plot, the absolute values of these differences
are found to be of the order of $10^{-9}$, so indeed the differences
between the two sets of power curves are negligible. The
practical value of the
equivalent short fiber calculation lies in the fact it gave 
essentially the same power curves about 100 times faster than the real-length fiber
calculation of Figure~\ref{fig:Tm_Yb_10m}.

In the second scenario, the total input power of 30~W is distributed
equally between the LP01 and LP11 modes. From the top right panel of
Figure~\ref{fig:Tm-individual}, we find that LP01 mode amplifies more
than the LP11 mode. Moreover, as in the left panel, the results from
the real and equivalent short fiber are visually
indistinguishable. However, a more careful examination of the
difference $P_l \circ \zeta - \tilde{P}_l$ in the bottom right plot
shows that maximal absolute power differences are about 0.3 near the
inlet of the fiber. Although this is many fold larger than the first
scenario, the relative power error of $3 \times 10^{-4}$ is still
quite small enough to make the equivalent short fiber a useful
practical tool. Note that the difference
$P_l \circ \zeta - \tilde{P}_l$ is now highly oscillatory, due to the
interactions between the two modes.

\subsection{Realizing the equivalent short fiber for the Yb-doped case}
\label{ssec:Yb-equiv-realize}

The equivalent short fiber in the Yb-doped case is more easily
realizable than the Tm-case as the Yb population dynamics is simpler.
The following conclusion can be arrived at easily proceeding similarly
as in Subsection~\ref{ssec:Tm-equiv-realize}.
\begin{center}
  \framebox{\parbox{\textwidth}{ %
      A short fiber of length $\tilde L$ is equivalent to a Yb-doped
      fiber of length $L$ if the fiber's original parameter $\Nt$ is
      changed to $\tnt = L \Nt / \tilde L$, i.e., this change
      realizes~\eqref{eqn:equiv cmt}.  }}
\end{center}

\begin{figure}
  \centering
  \includegraphics[width=0.5\textwidth]{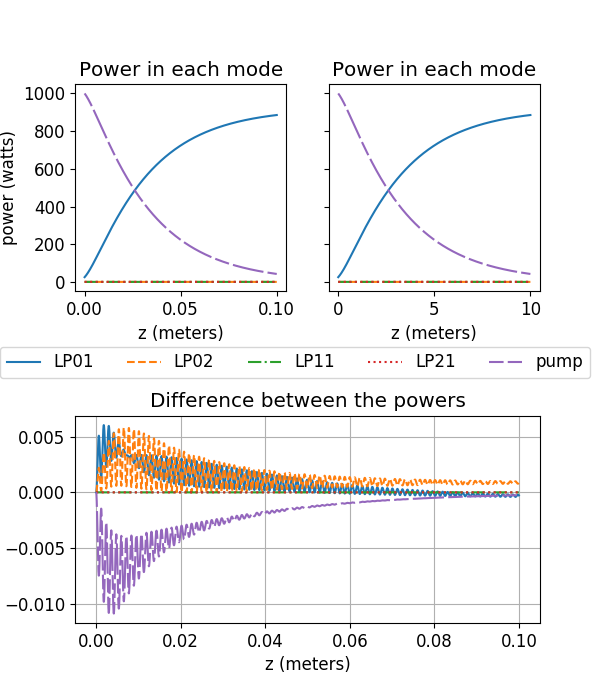}
  \hspace{-0.2cm}
  \includegraphics[width=0.5\textwidth]{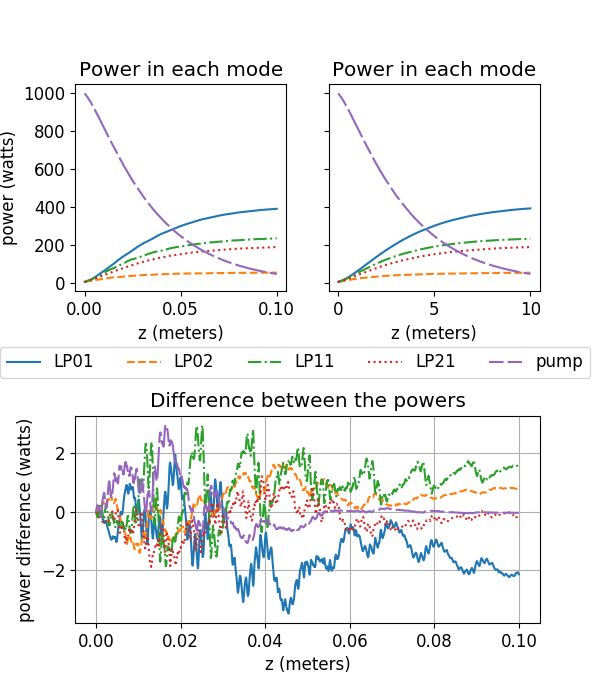}
  \caption{A comparison between a {\bf{Yb-doped fiber}} and its
    equivalent short counterpart. The left panel shows the case where
    the input signal power was wholly contained in the LP01 mode,
    while the right panel shows the case where it was equally
    distributed between all four modes.}
  \label{fig:Yb-individual}
\end{figure}

Figure~\ref{fig:Yb-individual} gives some indication of the practical
performance of this equivalent short fiber. As in the experiments for
the Tm-fiber reported in Figure~\ref{fig:Tm-individual}, here we
consider two scenarios, the first where all input signal power is
given to the LP01 mode, and the second where the input power is
distributed to the four LP modes equally (25\% each). The left panel
in Figure~\ref{fig:Yb-individual} shows the former, while the right
panel shows the latter. The equivalent fiber is less faithful in the
latter case, but the scale of the errors observed in the bottom plots
in both cases are well within the acceptable error ranges in
engineering practice. {(Laboratory power measurement uncertainties tend to be about $\pm5\%$.)}

\subsection{Increase of error with respect to some parameters}

We want to understand how relative power differences between the
equivalent and real fiber vary with respect to two important input
parameters $P_p^0$ and the short fiber length $\tilde L$. We consider
both the Tm and Yb fibers, holding the
original fiber length $L$ fixed to 10~m.

The solutions of the original and equivalent fiber models vary as
initial conditions are changed. Therefore to compare one with the
other in the {\em worst case} scenario, we take the maximum of the
power error measures over the set
\[
  \cA = \left\{ \alpha \in \mathbb{C}^M:  \; \int_{\om_z}
    I_s(x, y, 0, \alpha) \; dx dy = P_s^0 \right\},
\]
i.e., the set $\cA$ is the set of all input distributions yielding the
same initial signal power~$P_s^0$, which is set for Tm and Yb fiber
per Tables~\ref{tab:Tm} and~\ref{tab:Yb}, respectively. The initial
pump power $P_p^0$ is varied in the range 1000--5000~W (thus providing a
corresponding range of initial values for the $I_p$-component in the
model).  We solve the full CMT model and the equivalent short fiber
model, not only for this range of $P_p^0$, but also for decreasing
values of the short fiber length $\tilde L$.  The following quantity
is then computed across all such solutions:
\begin{equation}
  \label{eq:eps-def}
  \varepsilon(P_p^0, \tilde{L}) =
  \max_{A^0 \in \cA }\; 
  \frac 
  {\displaystyle \max_{l=0,1, \ldots, M} \; \max_{0 \le z \le L}
  \big| (P_l -  \tilde{P}_l \circ \zeta^{-1} )(z) \big|}
  {\displaystyle \max_{l=0,1, \ldots, M} \;\max_{0 \le z \le L}
    \big| P_l(z)\big| }.  
\end{equation}
Thus $\varepsilon$ represents the maximal possible power deviations
between the equivalent and original models over all input signal
distributions and over all mode components, as a function of initial
pump power $P_p^0$ and the fictitious length $\tilde{L}$.  Values of
$\varepsilon$ will thus inform us of the ranges of $P_p^0$ and
$\tilde L$ where the equivalent short fiber is more useful.

To practically compute $\veps$, we replace the maximum over the
infinite set $\cA$ by a computable maximum over a finite set obtained
by assigning each mode component all possible values from 0 to 100\%
in 10\% increments (while constraining the total signal power to
$P_s^0$).  In the case of the 2-mode thulium fiber, this resulted in
11 input power distributions, %
while for the ytterbium-doped fiber having 4 modes, 286 distributions
were required.  The maximum over $z$ in~\eqref{eq:eps-def} is replaced
by the maximum over the points where ODE solver traversed.  We used
polynomial degree $p=5$ for the finite element approximation of modes
and the 7-stage Dormand-Prince Runge Kutta method for solving the ODE
system. Collecting data from hundreds of simulations, we then plot
$\veps$ in a two-dimensional grid of $P_p^0$ and $\tilde L$ values.

\begin{figure}
  \centering
    \includegraphics[scale=.39]{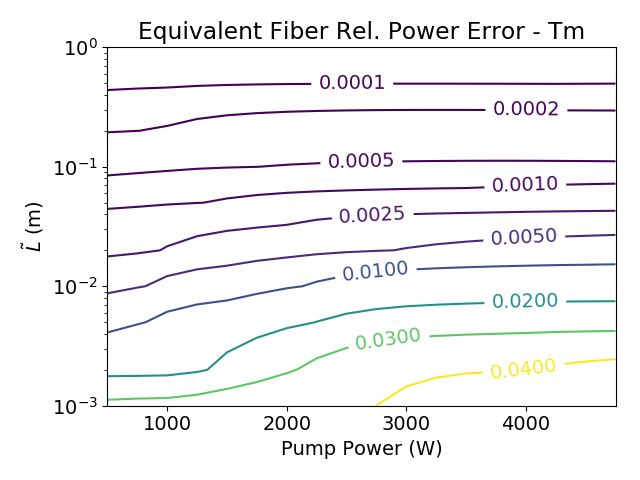}\hspace{-0.3cm}
    \includegraphics[scale=.39]{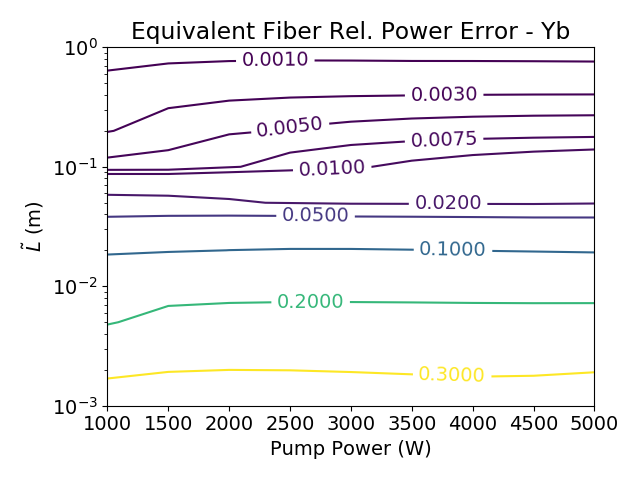}  
      \caption{
        Maximal relative power differences between a 10~m long real
        fiber and equivalent short fibers of various lengths $\tilde L$,
        for various initial pump powers $P_p^0$. The Tm case is shown on
        the left and Yb case on the right.
      }
  \label{fig:rel_pwr_diff}
\end{figure}

The resulting contour plots of the function $\varepsilon$ are given in
Figure~\ref{fig:rel_pwr_diff} for Yb and Tm fibers, for a range of
$P_p^0$ and $\tilde L$ values. We find that relative error
$\varepsilon$ varies mildly with respect to $P_p^0$ for any fixed
$\tilde L$, indicating that the absolute error in the powers increases
more or less linearly as $P_p^0$ is increased.  Looking vertically at
the plots of Figure~\ref{fig:rel_pwr_diff}, we find that holding
$P_p^0$ fixed, there are significant variations in $\varepsilon$ with
respect to $\tilde L$. The errors definitively increase as $\tilde{L}$
decrease. Figure~\ref{fig:rel_pwr_diff} clearly indicates that
excessively short equivalent fiber lengths are not recommendable.


\end{document}